\documentclass[12pt,a4paper]{amsart}
\usepackage{euscript,amsfonts,amssymb,amsmath,amscd}

\usepackage[T2A]{fontenc}

\usepackage[utf8]{inputenc}        

\usepackage[english,russian]{babel}

\usepackage{hyperref}

\input{epsf}

\newcommand{\m}{\mathbb}

\theoremstyle{plain}
\newtheorem{theorem}{Theorem}
\newtheorem{lemma}{Lemma}

\theoremstyle{definition}

\begin{document}

\title {Almost complex circle actions with few fixed points}

\author{Andrey Kustarev}

\renewcommand{\abstractname}{\empty}

\begin{abstract}We show that almost complex circle actions with exactly three fixed points do not exist in dimension $8$ and present an infinite series of $6$-dimensional manifolds possessing an almost complex circle action with exactly two fixed points.
\end{abstract}

\maketitle

Let $M$ be a smooth manifold with a smooth action of Lie group.
Since Conner and Floyd (\cite{kf}) theory of smooth actions of Lie groups has extensively involved the use of extraordinary cohomology theories, formal group laws and localization theorems for both ordinary and equivariant cohomology. The subject of this paper are circle actions preserving some almost complex structure on given manifold and having only isolated fixed points. Recently there has been progress on this and related problems (\cite{pelayotolman},\cite{tw},\cite{symp3pt}).

As shown in \cite{pelayotolman}, if $M^{n}$ is a manifold with an almost complex circle action and odd number of isolated fixed points, then $n=4k$. The case of $M^4=\m CP^2$ provides us with an example of smooth almost complex action with exactly three fixed points. 

\begin{theorem}
\label{T1}
There does not exist 8-dimensional manifold $M$ with a smooth almost complex circle action having exactly three isolated fixed points.
\end{theorem}

Recall that quaternionic projective space $\m HP^2$ admits an action of circle with exactly three fixed points but does not carry any almost complex structure. 

If smooth action of $S^1$ on $n$-dimensional manifold is almost complex and has exaclty two fixed points, then $n=2$ or $6$; furthermore, weights in fixed points of the action are $\{a,b,-a-b\}$ и $\{-a,-b,a+b\}$ for some positive integers $a,b$ if $n=6$. The proof may be found in \cite{pelayotolman}. The action is assumed to be symplectic in original paper but is extended easily to almost complex case. The main tool is the localization theorem for Chern numbers (\cite{ab},  see also \cite{toricgenera} for localization of toric genus).

\begin{theorem}
\label{T2}
For every almost complex manifold $M^4$ the manifold $(S^6)\,\#_{S^1}\,(S^1\times S^1\times M^4)$ admits an almost complex structure and, furthermore, circle action preserving that structure and having exactly two fixed points.
\end{theorem}

Here by $(S^6)\,\#_{S^1}\,(S^1\times S^1\times M^4)$ we mean fiberwise connected sum of manifolds along the tubilar neighborhood of two free orbits of $S^1$ action in each manifold.

{\it Proof of theorem \ref{T1}.} Let $M=M^8$ be a smooth manifold with a smooth almost complex circle action such that its fixed points are all isolated and there are exactly three of them; $a_1=m$ -- the maximum in absolute value among all weights of the action, $x\in M$ -- the fixed point corresponding to $a_1$. We may assume that $m>0$, otherwise switch to conjugate action.

Any $2k$-manifold with smooth almost complex circle action and weights equal to $\pm1$ has exactly $2^k$ fixed points if they are all isolated (\cite{tw},\cite{feldman}). Hence $m\ne\pm1$. Furthermore, the weights of $S^1/_{\m Z/m}$ acting on manifold $M^{\m Z/m}$ are also $\pm1$ so the number of fixed points of $S^1/{\m Z_m}$-action is also the degree of two. 

The number of fixed points may not be equal to one -- it would contradict with localization theorem for Chern numbers. This implies that component of $M^{\m Z/m}$ containing fixed points of the action is two-dimensional.  Let $y\ne x$, $y\in M^{\m Z/m}$ be a second fixed point such that  $b_1=-a_1=-m$ ( $b_i$ are weights in $y$). Then sets of weights $(a_2,a_3,a_4)$ and $(b_2,b_3,b_4)$ are equal modulo $m$ since the residuals are weights of $\m Z/m\subset S^1$ acting in space of normal bundle of $M^{\m Z/m}$ component in $M$. We denote weights in third fixed point $z\in M$ by $c_i$.

\begin{lemma}
\label{L1}
\cite{musin}
The weights $a_1,\dots,c_4$ may be split in pairs of weights opposite in sign and equal by absolute value; moreover, weights in any pair may be assumed to belong to different points.
\end{lemma}

\begin{lemma}
\label{L2}
\cite{pelayotolman}
If the number of fixed points of almost complex circle action on $M^{2n}$ is less than $n+1$, then for every fixed point $x\in M^{2n}$ exists another fixed point having same sum of weights as $x$.
\end{lemma}

Lemma \ref{L2} implies that sums of weights in all fixed points are equal. Lemma \ref{L1} implies that sum of all weights is zero and therefore the sum in every distinct point is zero. This means that weights $b_2,b_2,b_4$ are of the form $a_2+m, a_3+m, a_4$. The weight $a_4$ may not be opposite to any of weights $a_2+m, a_3+m$ because in that case we would have $m+a_2+a_3+a_4\ne 0$. Moreover, $a_2+a_4\ne 0, a_3+a_4\ne 0$. So the only remaining opportunity is $c_1=c_2=-a_4$.
Using similar arguments, we conclude that $a_2+m=-a_2, c_3 = -a_3, c_4 = -a_3-m$. 
But then manifold $M^{\m Z/(m/2)}$ (more exactly, its component carrying fixed points of the action) is four-dimensional and has exactly two fixed points. The weights of  $S^1/\m Z/(m/2)$ action on this manifold are $(2,-1)$ and $(-2,1)$ -- this contradicts localization theorem. $\Box$

{\bf Remark.} Here is an example of three sets of <<weights>> -- $$(-7,-1,10,12), (-8,3,5,14), (-8,3,5,14)$$ -- that satisfy every condition of localization theorem but still can't be weights of some almost complex circle action on $M^8$. Note that sum of weights is equal to $14$ in every set, so the implication of lemma \ref{L2} still holds and statement of lemma \ref{L1} is obviously false.

{\it Proof of theorem \ref{T2}.} It is well-known that $S^6 = G_2/SU(3)$ admits an almost complex structure that is preserved under action of maximal torus $T^2\subset SU(3)$ (\cite{borelhirz}). There are exactly two fixed points and one may find one-dimensional abelian subgroup acting on $S^6$ with weights $\{a+b,-a,-b\}$ and $\{-a-b,a,b\}$. 

Let $M^6$ and $N^6$ be manifolds with almost complex circle actions, $U_M$ and $U_N$ -- tubilar neighborhoods of some free orbits of actions in $M^6$ and $N^6$. Fix two equivariant diffeomorphisms $\phi_M, \phi_N : S^1\times S^4 \to \partial U_M, \partial U_N$ and denote by $X$ the smooth manifold $(M^6-U_M)\cup_{\phi_M} (S^1\times S^4\times I) \cup_{\phi_N} (N^6-U_N)$.
Theorem \ref{T2} is therefore reduced to the following statement.

\begin{lemma}
\label{L6}
The manifold $X$ admits almost complex structure that is itself an extension of structures on 
 $M^6-U_M$ and $N^6-U_N$. This structure is invariant under natural extension of $S^1$-action. 
\end{lemma}

{\it Proof of the lemma.} Theorem on existence of equivariant tubular neighborhood (\cite{bredon}) allows to extend smooth structure from $M^6-U_M$ and $N^6-U_N$ to $X$. 
The obstruction to equivariant extension of almost complex structure lies in the group $\pi_4(SO(6)/U(3)) = \pi_4(SO/U)$. This group is trivial (\cite{bott}).

The author is grateful to V. M. Buchstaber, T. E. Panov and also to A. Gaifullin, E. Grechnikov and Yu. Ustinovsky for valuable discussions.

\renewcommand{\refname}{References}

\end{document}